\documentclass{article}

\usepackage{a4}

\usepackage[utf8]{inputenc}

\usepackage{epsfig}

\begin{document}

\title{Counterexamples in non-positive curvature}

\author{Yves Coudene, Barbara Schapira}

\date{09/04/2010}

\maketitle

\begin{center}

\emph{Université de Bretagne Occidentale, 6 av. Le Gorgeu, 
29238 Brest cedex, France}

\emph{LAMFA, Universit\'e Picardie Jules Verne, 33 rue St Leu, 
80000 Amiens, France}
\end{center}

\begin{abstract}
  We give examples of rank one compact surfaces on which there exist recurrent
  geodesics that cannot be shadowed by periodic geodesics. We build rank one
  compact surfaces such that ergodic measures on the unit tangent bundle of the
  surface are not dense in the set of probability measures invariant by the
  geodesic flow. Finally, we give examples of complete rank one surfaces for
  which the non wandering set of the geodesic flow is connected, the periodic
  orbits are dense in that set, yet the geodesic flow is not transitive in
  restriction to its non wandering set. \footnote{37B10, 37D40, 34C28}{}
\end{abstract}

\section{Introduction}

Geodesic flows on negatively curved manifolds have been extensively studied
both from the topological and the measurable point of view.
Let $M$ be a complete connected pinched negatively curved riemannian manifold,
and denote by $g_t: T^1M\rightarrow T^1M$ the geodesic flow defined 
on the unit tangent bundle $T^1M$ of $M$. The non-wandering set of $g_t$
is denoted by $\Omega$.
$$
\Omega=\{\, v\in T^1M \ |\ \forall \, V \hbox{ neighborhood of } v, \ 
\exists \, t_n\rightarrow +\infty 
\hbox{ such that } g_{t_n}(V)\cap V \neq \emptyset \,\}
$$
Recall from the Poincaré recurrence theorem that $\Omega=T^1M$ if $M$ is compact
or with finite volume.

\medskip

\quad
Here are some basic properties of $g_t$.

\medskip

\emph{$\bullet$ 
Every vector in $\Omega$ can be shadowed by a
periodic geodesic. As a consequence, periodic geodesics are dense in the 
non wandering set of the flow.}

\medskip

\emph{$\bullet$ 
Assuming $\Omega= T^1M$, there are ergodic measures with full support that
are invariant by $g_t$. In fact ergodicity and full support are generic
properties in the set of all invariant probability measures.}

\medskip

\emph{$\bullet$ 
The flow $g_t$ is transitive in restriction to $\Omega$
 as soon as $\Omega$ is connected.}

\medskip

\quad
There has been several attempts to generalize these results to the 
case of geodesic flows defined on non-positively curved manifolds.
The goal of this paper is to provide explicit counterexamples to each
of the results stated above, in the context of rank one manifolds.

\medskip

\quad 
Let us recall what is a rank one manifold. A vector $v\in T^1M$ is a
\emph{rank one vector} if the only parallel Jacobi fields along the geodesic
generated by $v$ are proportional to the generator of the geodesic flow. A
connected complete non-positively curved manifold is said to be a 
\emph{rank one manifold} if its tangent bundle admits a rank one vector. 
In that case, the set of rank one vectors is an open subset of $T^1M$.

\medskip

\quad
A vector $v$ on a non-positively curved surface has rank one if and only if
there is a point of negative sectional curvature on the geodesic generated by
$v$. Hence, a rank one surface is simply a non-positively curved surface with at
least a point where the curvature is negative.

\medskip

\quad
We shall build rank one surfaces $M$ such that

\medskip 

\emph{$\bullet$ the surface $M$ is compact, and there exist recurrent
  rank one geodesics that cannot be shadowed by periodic geodesics;}

\medskip

\emph{$\bullet$ the surface $M$ is compact, and ergodic measures are not dense
  in the set of all probability measures invariant by $g_t$;}

\medskip

\emph{$\bullet$ the set $\Omega$ is connected, the periodic orbits are dense
  in $\Omega$, and the flow $g_t$ is not transitive in restriction to $\Omega$.}

\medskip

\quad
All these counterexamples contain embedded flat cylinders.
It is plausible that in dimension two, such cylinders are the only obstruction
to the results stated above. In higher dimension however, the question 
is wide open.

\medskip

The article is organized in three sections, each devoted to one of the
counterexamples stated above.

\section{The closing lemma}

The first instance of a closing lemma appears in the work of J. Hadamard 
\cite{had} on negatively curved surfaces embedded in ${\bf R}^3$.
J. Hadamard showed that a piece of geodesic coming back close to its 
starting point is in fact shadowed by a periodic closed geodesic.
As a consequence, the set of periodic closed geodesics is dense 
in the closure of the set of recurrent geodesics. 

\medskip

\quad 
If the negatively curved surface is of finite volume, then the celebrated
theorem of Poincaré ensures that the recurrent geodesics are in fact dense in
the manifold. So periodic vectors are dense in the unit tangent bundle of the
surface. Remarkably, the result of J. Hadamard is contemporary to the Poincaré
recurrence theorem, and there were no examples of finite volume negatively
curved surfaces embedded in ${\bf R}^3$ at that time. First examples of such
surfaces are attributed to Vaigant in \cite{geom}, although we were unable to
locate the original reference. Yet J. Hadamard provided examples of negatively
curved embedded surfaces in ${\bf R}^3$, and the closing lemma was one of
the key ingredient in the study of recurrent geodesics on these surfaces.

\medskip

\quad Purely dynamical proofs of the closing lemma were given by D.V. Anosov, who
extended that result to the class of systems that now bear his name. The closing
lemma has become a classical tool in the study of systems exhibiting some form of
hyperbolic behavior, and the name of Anosov is frequently associated with the
theorem. In the context of negatively curved manifolds, the Anosov closing lemma
reads as follows.

\medskip

{\bf Theorem 1 } (Anosov closing lemma \cite{an})

\emph{Let $M$ be a compact negatively curved manifold, $g_t$ the geodesic flow
on $T^1M$.
Then for all $\varepsilon>0$, there exists $\delta>0$
and $T_0>0$, such that for all $w\in T^1M$, we have:
\vspace{.5em}\\
for all $T\geq T_0$ such that $d\bigl(g_T(w),(w)\bigr)<\delta$,
there exists a periodic vector $w_0$
of period $l$ such that
$$
|T-l|<\varepsilon
$$
$$
\forall t\in [0,T],\quad d\bigl(g_t(w),g_t(w_0)\bigr)<\varepsilon.
$$
}

\medskip

\quad
The closing lemma implies the density of the periodic orbits in the
non-wandering set of the flow. The term ``Anosov closing lemma'' is sometimes
used to describe that last property \cite{rob}. We will see below examples of
(non-uniform) hyperbolic systems for which the closing lemma is not satisfied,
yet the periodic orbits are dense in the ambient space.

\medskip

\quad
Note also that the compactness hypothesis on $M$ allows for uniform
$\delta$, $\varepsilon$. We will need some local version of the closing lemma.

\bigskip

{\bf Definition}

\emph{We say that the closing lemma is \emph{satisfied around a vector 
$v\in T^1M$} if there is a neighborhood $V$ of $v$
such that:}

\medskip

\emph{$\forall \, \varepsilon >0, \ \exists \, \delta>0,\ \exists\, T_0>0, 
\ \forall\, w\in V,\ \forall\, T\geq T_0,$}

\smallskip

\emph{the conditions $g_T(w)\in V$ and $d\bigl(g_T(w),w\bigr)<\delta$ 
imply that there is a periodic vector $w_0$ of period $l$ satisfying 
$|T-l|<\varepsilon$ and $d\bigl(g_t(w),g_t(w_0)\bigr)<\varepsilon$ 
for all $t\in [0,T]$.}

\bigskip

\quad We now turn to non-positive curvature. A version of the closing lemma is
given by W. Ballmann, M. Brin, R. Spatzier in \cite{bbs}, th. 4.5. The authors
are mainly interested in higher rank manifolds, although their arguments hold in
the rank one setting. They show that the closing lemma is satisfied around every
vectors of minimal rank, and the proof carries to the rank one setting. 
As an example, the reader may want to convince himself by an
elementary argument that the closing lemma indeed holds on a flat torus. In that
example, all vectors have a rank equal to the dimension of the torus.

\medskip

\quad
Let us now restrict the discussion to rank one manifolds.
The following result is stated in
\cite{eb1}, th. 4.5.15. We refer also to \cite{cs} for a shorter proof on rank
one manifolds.

\bigskip

{\bf Theorem 2 } (Rank one closing lemma \cite{eb1})

\emph{Let $\tilde{M}$ be a rank one manifold.
Then the closing lemma is satisfied around all rank one vectors.}

\bigskip

\quad This result implies that on a compact rank one manifold, periodic vectors
are dense in the set of rank one vectors. That set is both open and dense, hence
we see that periodic vectors are in fact dense in $T^1M$. This can be extended
to rank one manifolds for which $\Omega = T^1M$.

\medskip

\quad 
Establishing the density of periodic orbits is the first step in the proof of the
transitivity of the geodesic flow on compact connected rank one manifolds. We
will use that property in the sequel, so let us state the following result due
to P. Eberlein \cite{eb}.

\bigskip

{\bf Theorem 3 } \cite{eb}

\emph{Let $M$ be a compact rank one manifold.
Then the geodesic flow on $T^1M$ is transitive: there exists a vector in $T^1M$
with a dense orbit.}

\bigskip

\quad
Without the compactness assumption on $M$, transitivity holds as
soon as $\Omega=T^1M$  \cite{eb1} (1.9.15, 4.7.4,
4.7.3).

\bigskip 

Let us now turn our attention to compact non-positively curved surfaces. We have
seen that the closing lemma is satisfied on flat toruses. Also, on all other
surfaces, there is an open dense set of vectors around which it is satisfied.
Yet there are examples of surfaces on which there exists recurrent rank one
trajectories that cannot be closed. In fact, that happens as soon as there is an
embedded flat cylinder in the surface.

\medskip

We say that a surface $M$ contains an \emph{embedded flat cylinder} if there is
a compact subset of $M$ which is isometric to the euclidean cylinder
$[0,l]\times {\bf R}/(2\pi h{\bf Z})$, for some $l,h>0$. Geodesics in the
cylinder are straight lines. The periodic geodesics contained in the cylinder
are parallel among themselves, and define a direction in the tangent bundle of
the cylinder, that will be deemed as ``vertical''. This direction will be
depicted vertically in the figures.

\medskip

\quad
Note that all negatively curved surfaces carrying a simple closed geodesic can
be flattened around that geodesic so as to obtain a surface with an embedded
cylinder. 

\medskip 

\quad
We can now state the counterexample to the closing lemma. 

\bigskip

{\bf Theorem 4 }
 
\emph{Let $M$ be a complete surface with an embedded flat cylinder. We assume
  that the geodesic flow is transitive on $T^1M$. Then the closing lemma is not
  satisfied around vectors generating closed geodesics contained inside the
  cylinder.}

\bigskip

\quad
No hypothesis on the curvature of $M$ is needed.
The transitivity assumption rules out the case of a flat torus, where
the closing lemma holds everywhere.
Transitivity is satisfied in the rank one setting, so we get:

\bigskip

{\bf Corollary}
 
\emph{The Anosov closing lemma does not hold on compact rank one surfaces
  admitting an embedded cylinder.}

\bigskip

{\bf Proof}\\
Let $v$ be a periodic vector whose trajectory is contained in the interior
of the cylinder. The trajectory of $v$ bounds two connected components of 
the cylinder, denoted by $C_1$ and $C_2$. 

\medskip 

\quad Given $\varepsilon>0$ and $\theta\in ]0,\pi/2[$, we consider the open set
$U_1$ consisting in vectors whose base points are in $C_1$, at distance less
than $\varepsilon$ from the base point of $v$, and whose angle with the vertical
direction belongs to $]-\theta,0[$. We also consider the open set $U_2$ of
vectors with base points in $C_2$ close to $v$, and with angles belonging to
$]0,\theta[$.

\bigskip

\centerline{
\epsfig{figure=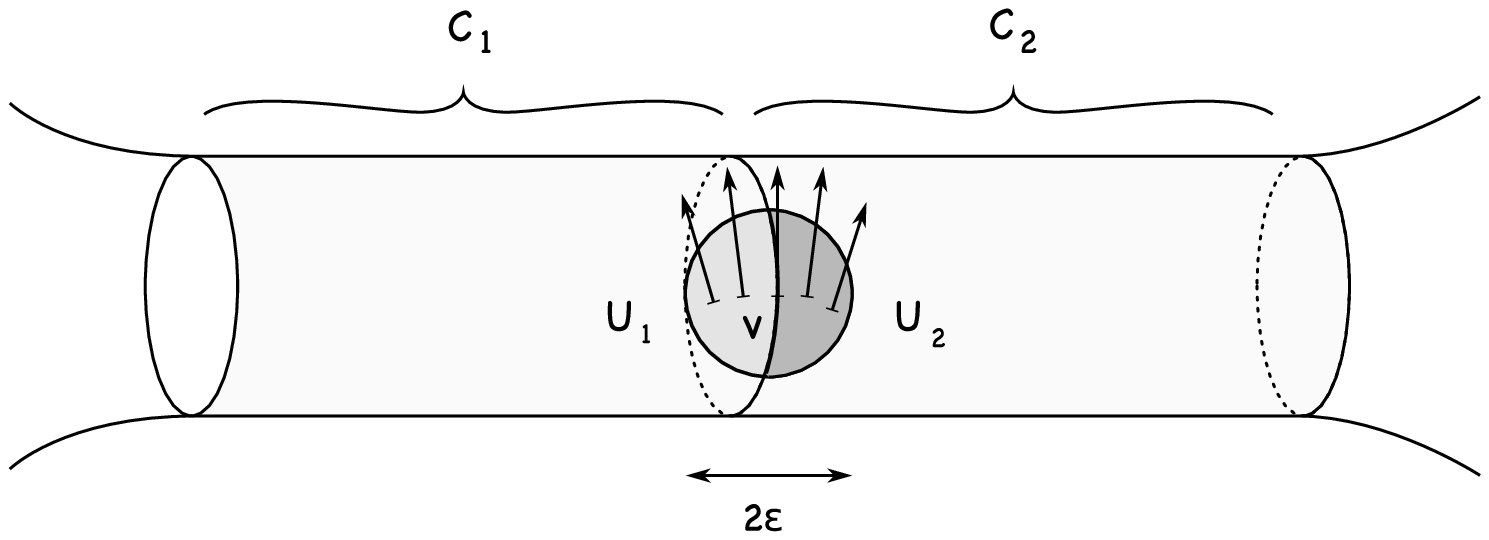,width=.8\textwidth,angle=0}}

\bigskip

\quad Recall that the geodesic flow is transitive. There is a vector $w\in T^1M$
with a dense orbit. The orbit of $w$ crosses the closed geodesic generated by
$v$ infinitely many times. The geodesic generated by $w$ enters $U_2$ at time
$t_0$ and then enters $U_1$ infinitely often. This means that we can find times
$t_n>0$ arbitrarily large, so that $d(g_{t_0}w,v)$ is less than
$\theta+\varepsilon$ and $d(g_{t_n}w,g_{t_0}w)$ is less than
$2\theta+2\varepsilon$. See next figure.

\bigskip

\centerline{
\epsfig{figure=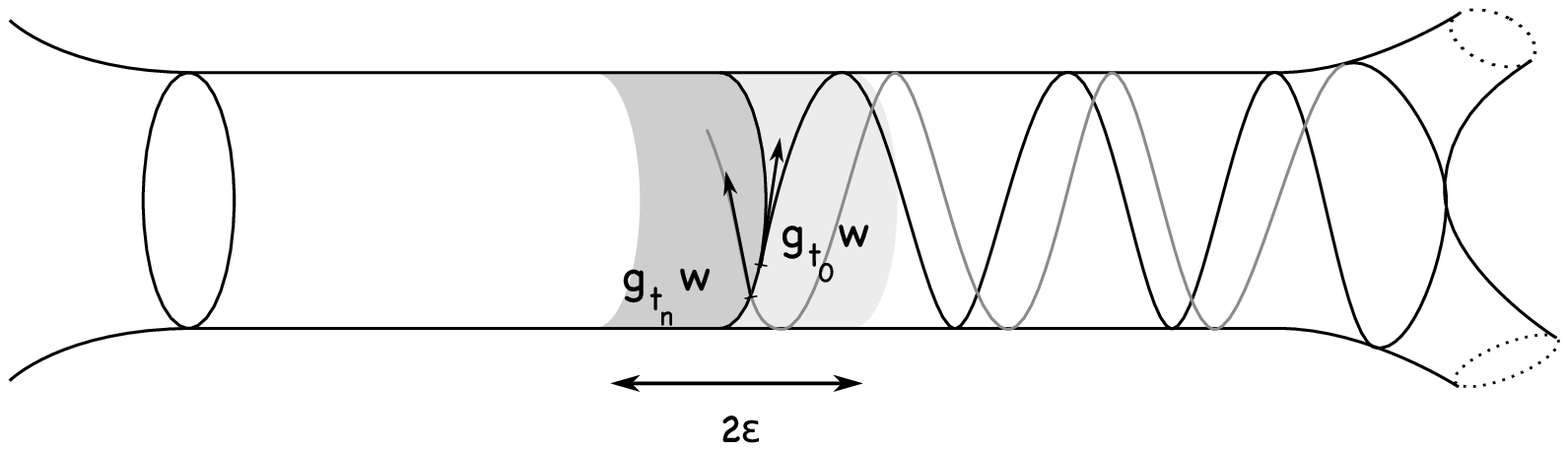,width=.8\textwidth,angle=0}}

\bigskip

\bigskip

Orbits in the cylinder are straight lines. Hence an orbit of the geodesic flow
entering the cylinder by one side must leave the cylinder by the other side.

\medskip

\quad
If the closing lemma is true in the neighborhood of $v$, we can find a periodic
orbit $w_0$ with period close to $t_n$,
that starts close to $g_{t_0}(w)$ and follows the orbit of $g_{t_0}(w)$ until
time $t_n$. 
The orbit $g_t(w_0)$ must leave the cylinder by the right side for $t>0$,
since it follows the orbit of $g_{t_0}(w)$ for positive $t$.
It must also leave the cylinder by the right side for $t<0$, 
since it follows the orbit of $g_{t_n}(w)$ for negative $t$.
This is in contradiction with a closing lemma around $v$.

\bigskip

\bigskip

\centerline{
\epsfig{figure=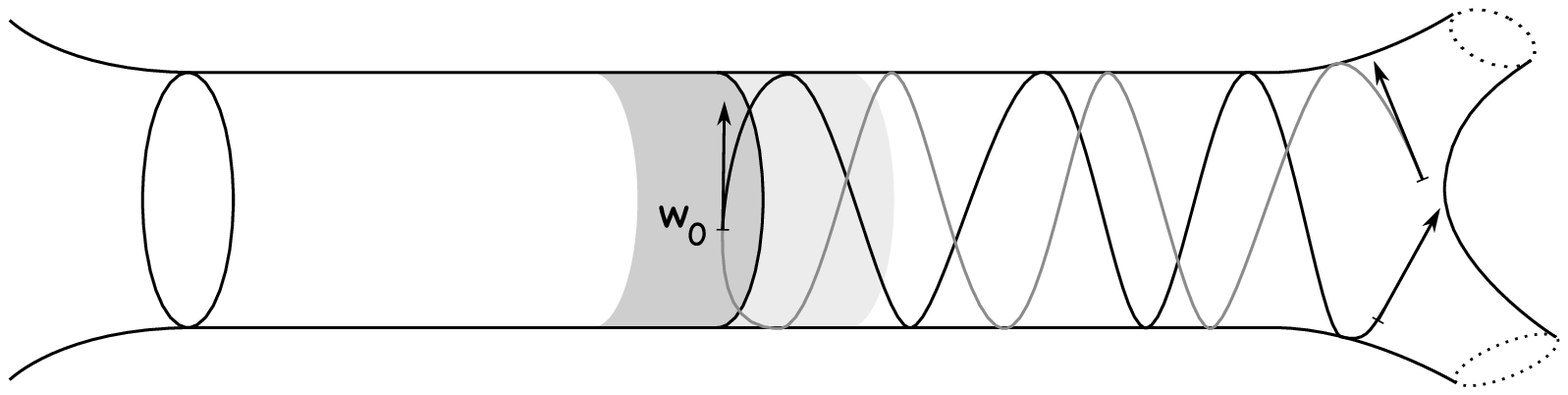,width=.8\textwidth,angle=0}}

\bigskip

This ends the proof of the result. 

\medskip

\quad The same argument would show that there is no local product structure in
the neighbourhood of $v$: it is impossible to find an orbit close to $v$ and
$g^{t_0}w$, negatively asymptotic to the orbit of $g^{t_0}w$ and positively to
the orbit of $g^{t_n}w$ (see \cite{cs} for definitions).

\medskip

\quad
Finally, we remark that the closing
lemma holds around $g_{t_0}w$ (see \cite{cs}). The orbit of $g_{t_0}w$ itself
must come back sufficiently close to $g_{t_0}w$ however, 
so as to be directed toward the same boundary of the cylinder as $w$, 
for the closure to happen.

\section{Genericity}

The closing lemma is a key ingredient in the understanding of the 
generic properties of measures on $T^1M$ invariant by the geodesic flow.

\medskip

\quad
Let $X$ be a complete metric space. Recall that a property is \emph{generic}
if the set of points satisfying that property contains a countable intersection
of open dense sets. A countable intersection of generic sets is again generic.
From the Baire category theorem, a generic set must be dense.

\medskip 

\quad 
The first use of the Baire category theorem in the study of 
geodesics is again due to J. Hadamard \cite{had}.
This allowed him to show that on a negatively curved surface embedded in 
${\bf R}^3$ with many periodic geodesics, there are
geodesics that accumulate on infinitely many periodic geodesics.

\medskip

\quad The following theorem of K. Sigmund, proven in the context of transitive
Anosov flows, can be seen as some quantitative strengthening of J. Hadamard
result.

\bigskip

{\bf Theorem 5 } \cite{si}\\
\emph{Let $M$ be a connected compact negatively curved manifold. Then the Dirac
  masses on periodic orbits are dense in the set of all probability measures on
  $T^1M$ invariant by the geodesic flow. Moreover, ergodicity and full support
  are generic properties in the set of invariant probability measures.}

\bigskip

\quad Note that there are indeed many explicit examples of ergodic probability
measures of full support, when the negatively curved manifold is compact. The
volume, the measure of maximal entropy, more generally the equilibrium state
associated to a H\"older potential, are all Bernoulli measures of full support.

\medskip

\quad
 We extended K. Sigmund's result to the rank one setting
in \cite{cs}, freeing ourselves from any compactness assumption. We denote by
${\cal R}_1$ the open set of rank one vectors in $T^1M$, and 
$ \Omega \cap {\cal R}_1$ the set of non-wandering rank one vectors on $T^1M$.
Our theorem can be stated as follows.

\bigskip

{\bf Theorem 6 } \cite{cs}

\emph{Let $M$ be a connected complete non-positively curved manifold. We assume
  that the geodesic flow on $T^1M$ admits more than two rank one periodic
  orbits. Then the set of ergodic measures on $ \Omega \cap {\cal R}_1$, with
  full support in $ \Omega \cap {\cal R}_1$, is generic in the set of all
  probability measures defined on $ \Omega \cap {\cal R}_1$.}

\bigskip

\quad
Strictly speaking, we proved that theorem for a subset 
$\Omega_1\subset \Omega \cap {\cal R}_1$
that contains the set of recurrent rank one vectors. Let us explain
why it is also true for $ \Omega \cap {\cal R}_1$. Note that the set of
recurrent rank one vectors is dense in $ \Omega \cap {\cal R}_1$,
this follows from the rank one closing lemma.
From the Poincaré recurrence theorem, an invariant
probability measure must give full measure to the set of recurrent vectors. As a
consequence, an invariant measure supported by ${\cal R}_1$ gives full measure
to ${\Omega}_1$. The set of invariant measures on $ \Omega \cap {\cal R}_1$
and $\Omega_1$ thus can be identified.

\medskip

\quad The assumption that there are more than two rank one periodic orbits rules
out the case of a hyperbolic cylinder. On such a cylinder, the non-wandering set
is made of two opposite periodic geodesics; it is not connected and the flow is
not transitive on that set.

\bigskip 

{\bf Corollary} \cite{cs}

\emph{Under the assumptions of the previous theorem,
  there exists an ergodic probability measure invariant by the
  geodesic flow, whose support contains all rank one periodic orbits.}

\bigskip

\quad So if we assume that $M$ is a connected complete rank one manifold
satisfying $\Omega=T^1M$, then there is always an invariant probability measure
of full support on $T^1M$. The existence of such measure is non trivial even
in the case of a surface with constant negative curvature.

\medskip

\quad The question of whether ergodicity and full support are generic properties
in the set of all invariant measures on $T^1M$ was left open in \cite{cs}. One
of the difficulties was the lack of a closing lemma valid for all vectors in
$T^1M$. Now, refining on the counterexample of the previous section, we shall
show that ergodicity is not a generic property in general.

\bigskip

{\bf Theorem 7 }

\emph{Let $M$ be a compact riemannian surface with an embedded flat cylinder.
  Then the Dirac measures supported by the periodic geodesics in the interior of
  the cylinder are not in the closure of the invariant ergodic probability
  measures with full support.}

\bigskip

{\bf Proof}\\
Let $c: {\bf R} \rightarrow M$  be a periodic geodesic lying inside the 
cylinder. Let us consider a tubular neighborhood around $c$ which is isometric
to $]-3\, \varepsilon , 3\, \varepsilon [\times S^1$ and let 
$\theta \in \ ]0,{\pi\over 2}[$. We denote by $U_\varepsilon\subset T^1M$ 
the neighborhood of $c$ containing all unit vectors tangent to 
$]-\varepsilon,\varepsilon\, [\times S^1$ whose angle 
with the vertical direction belongs to $]-\theta,\theta[$.
The set $U_\varepsilon$ is depicted below.

\bigskip

\centerline{
\epsfig{figure=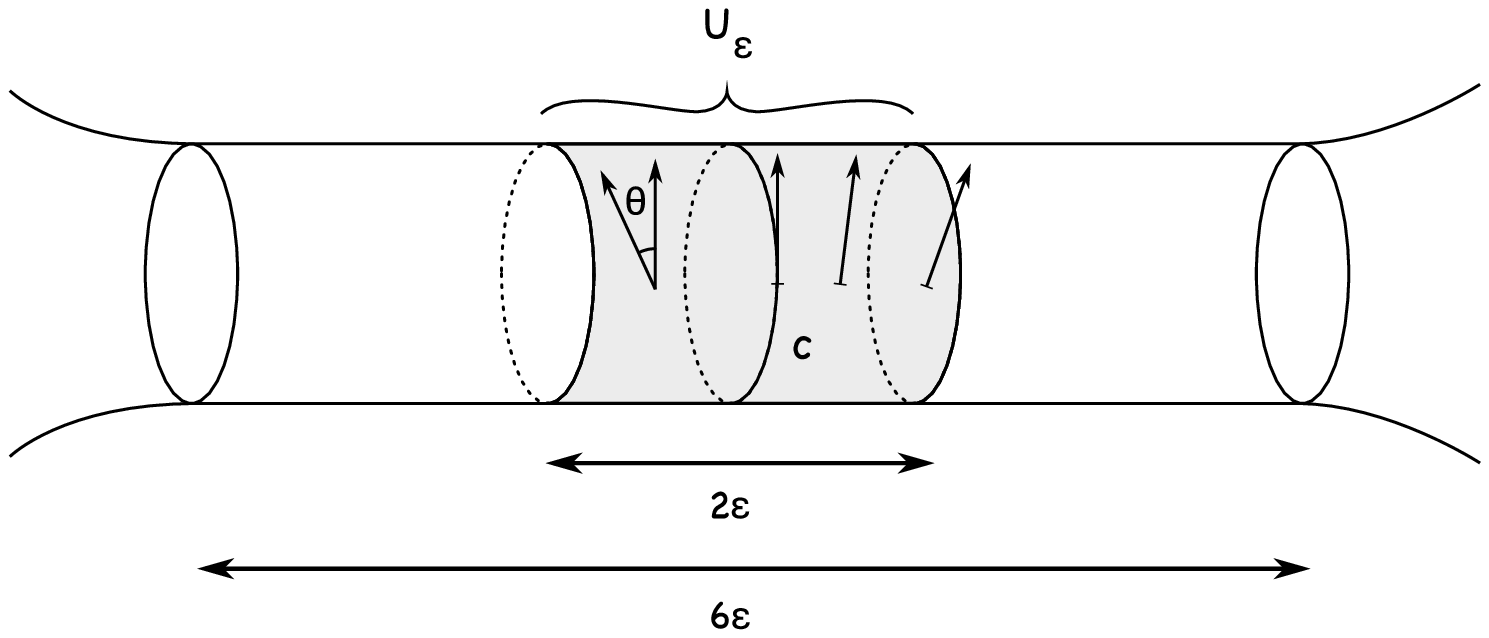,width=.8\textwidth,angle=0}}

\bigskip

The orbits of the geodesic flow in the cylinder are straight lines.
As a result, a trajectory that enters $U_\varepsilon$ at time $t_0$
and leaves $U_\varepsilon$ at time $t_1$, must have spent 
a time at least $t_1-t_0$ in the left cylinder 
$]-3 \varepsilon, - \varepsilon [\times S^1$ or in the right cylinder
$]\varepsilon, 3 \varepsilon [\times S^1$.
Hence, the quantity
$$
{1\over T} \ \lambda \Bigl(
\{\,t\in [0,T] \ |\ g_t(v)\in U_{\varepsilon}\, \}
\Bigr)
$$
is bounded from above by ${1\over 2}$ for all $T\geq 0$ and all $v\in T^1M$
that do not belong to the cylinder. See picture below.

\bigskip

\bigskip

\centerline{
\epsfig{figure=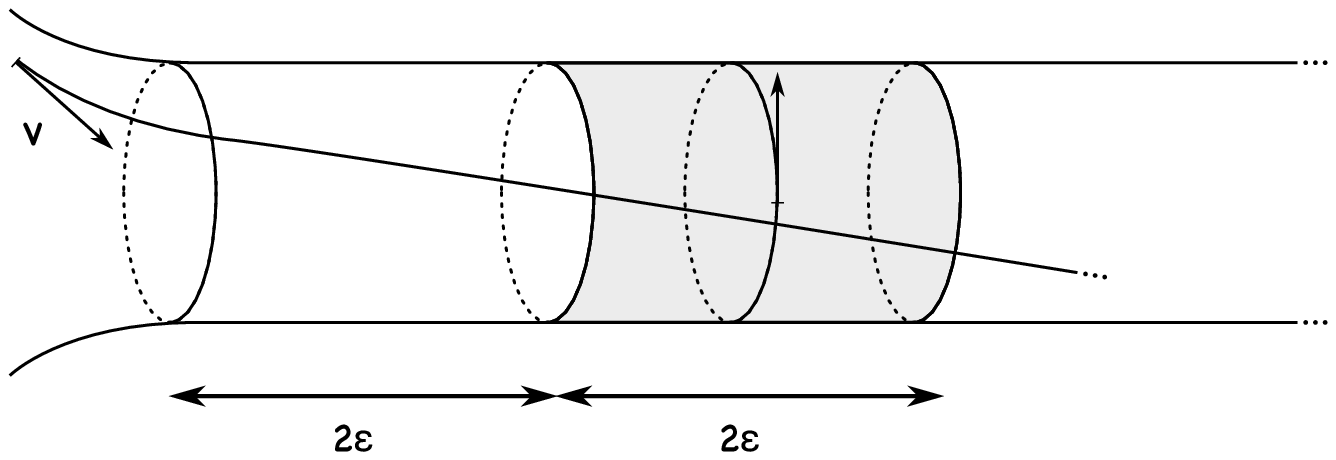,width=.8\textwidth,angle=0}}

\bigskip

Let $\mu$ be an ergodic measure with full support. If it is close enough
to the Dirac measure supported by $c$, then $\mu(U_\varepsilon)>{1\over 2}$.
Now we use the ergodic theorem. We can find some vector $v$ outside
the cylinder such that:
$$
{1\over T}\ \lambda \Bigl(
\{\,t\in [0,T] \ |\ g_t(v)\in U_{\varepsilon}\, \}
\Bigr) \ 
\mathop{\hbox to 4em{\rightarrowfill}}_{n\rightarrow\infty} \ 
\mu(U_\varepsilon)
$$
This gives a contradiction and the theorem is proven.

\bigskip

{\bf Remarks}


$\bullet$ Here again, there is no condition on the curvature of the
surface.

\smallskip

$\bullet$ The compactness assumption on $M$ plays no essential role in the proof.

\smallskip

$\bullet$ The same method can be applied in any dimension to build examples
of rank one manifolds on which ergodicity is not a generic property.

\bigskip

{\bf Corollary}\\
\emph{Let $M$ be a compact rank one surface with an embedded flat cylinder.
  Then ergodic measures are not dense in the set of all invariant probability
  measures on $T^1M$. Hence, ergodicity is not a generic property in $T^1M$.}

\bigskip

{\bf Proof}\\
Recall that the set of ergodic probability measures is a $G_\delta$ set. Also,
the set of measures with full support is a $G_\delta$-dense set if the periodic
orbits are dense in the ambiant space,
we refer to \cite{si} or \cite{cs} for a proof.

So, if the set of ergodic measures is dense, then the ergodic measures with full
support are dense by the Baire category theorem, and the (ergodic !) Dirac
measures in the cylinder can be approached by these measures, in contradiction
with the previous theorem. This proves the corollary.

\bigskip

\quad The reader may want an explicit example of a measure that is not in the
closure of the ergodic measures. A non trivial convex combination of Dirac
measures supported inside the cylinder does the job.

\medskip

\quad
Let us give a quantitative version of the previous theorem.
First, recall that the riemannian metric on $M$ induces a natural
riemannian metric on $TM$. This metric on $TM$ is uniquely caracterized
by the following properties:

\medskip

- the canonical projection from $TM$ to $M$ is a riemannian submersion,

\smallskip

- the metric induced in the fibers is euclidean,

\smallskip

- horizontal and vertical distributions are orthogonal.

\medskip

This metric induces in turn a metric on the submanifold $T^1M$. As an example,
the metric obtained on the unit tangent bundle $T^1C$ of a flat two dimensional
cylinder $C$ is just the euclidean metric on the product
$T^1C \simeq C\times S^1$.


\medskip 

The Prohorov metric on the set of Borel probability measures
is defined by
$$
\rho(\mu,\nu) = \inf\{\,\varepsilon>0 \ | \ 
\forall\,A \subset X,\ 
\mu(A)\leq \nu(V_\varepsilon (A))+ \varepsilon \, \}
$$

where $V_\varepsilon (A)$ is the $\varepsilon$-neighborhood of $A$,
$$
V_\varepsilon (A) = \{ y\in X \ | \ d(x,A) < \varepsilon \}.
$$
The metric $\rho$ is bounded by one.
We refer to \cite{bi} for its basic properties.

\bigskip 

{\bf Theorem 8 }

\emph{Let $M$ be a compact riemannian surface with an embedded euclidean
  cylinder isometric to $[0,l]\times {\bf R}/(2\pi h {\bf Z})$, $l,h>0$. Let
  $\mu$ be an ergodic probability measure of full support, invariant by the
  geodesic flow, $\delta$ the Dirac measure on a closed geodesic contained in
  the interior of the cylinder. The distance on $M$ between that closed geodesic
  and the boundary of the cylinder is denoted by $d$. Then
$$
\rho(\delta, \mu) \geq \min\Bigl(d, \, {l\over l+2}\Bigr)
$$}

\medskip 

This gives a bound of order $1-{2\over l}$ when $l$ is big.

\medskip

\quad
The bound does not depend on the height $h$ of the 
cylinder, nor on the geometry of $M$ outside the cylinder, nor on the 
measure $\mu$.
Here again we do not need any assumption on the curvature of $M$.

\medskip

\quad
For the geodesic in the middle of the cylinder, we have $d=l/2$ and the lower
bound is equal to $l\over l+2$.

\bigskip

{\bf Proof}\\
We take for $A$ the closed orbit in $T^1M$ supporting $\delta$.
Consider a vector $v$ in the cylinder, denote by $\theta$ the angle between
$v$ and the vertical direction, and $r$ the distance between the base point
of $v$ and the closed geodesic. The vector $v$ is in $V_\varepsilon (A)$
if and only if $\sqrt{\theta^2+r^2}$ is less than $\varepsilon$.

\medskip

\quad We now ensure that $V_\varepsilon (A)$ is contained in the cylinder by
assuming $\varepsilon < d$. This unfortunately rules out the case of a set $A$
given by a closed geodesic on the boundary of the flat cylinder. 
Without further assumption on the geometry of $M$, there is no way to guarantee
that $V_\varepsilon(A)$ does not contain arbitrarily long pieces of geodesics 
when $A$ is a bounding geodesic.

\medskip

\quad Let us consider a vector with a dense trajectory, that enters the cylinder
with an angle $\theta$ with respect to the vertical direction. It spends a time
$l/\sin(\theta)$ in the cylinder. If the trajectory crosses $V_\varepsilon(A)$,
that is if $\theta$ is less than $\varepsilon$, it spends a time
$2\sqrt{\varepsilon^2-\theta^2}/\sin(\theta)$ in $V_\varepsilon(A)$.
So the fraction of time spent in $V_\varepsilon(A)$ during the travel through
the cylinder is given by 
$$
2\sqrt{\varepsilon^2-\theta^2}\ \over l
$$
This can be bounded independently of $\theta$ by $2\varepsilon/l$, and this
bound is best when the vector enters the cylinder
close to the vertical direction.

\medskip 

\quad
Let $\mu$ an ergodic probability measure with full support. 
By the ergodic theorem, we can find some vector $v$ with dense trajectory
such that 
$$
{1\over T}\ \lambda \Bigl(
\{\,t\in [0,T] \ |\ g_t(v)\in V_{\varepsilon}(A)\, \}
\Bigr) \ 
\mathop{\hbox to 4em{\rightarrowfill}}_{n\rightarrow\infty} \ 
\mu(V_\varepsilon(A))
$$

\medskip

This gives $\mu(V_\varepsilon(A))\leq {2\,\varepsilon/l}$. Let us denote
by $\delta$ the Dirac measure on the periodic orbit under consideration.
Coming back to the definition of the Prohorov distance, the equation 
$$
1=\delta(A)\leq \mu(V_\varepsilon(A))+\varepsilon
$$
implies
$$
1\leq {2\over l} \ \varepsilon + \varepsilon 
$$
that is,
$$
{l\over l+2} \leq \varepsilon
$$
This estimate is greater than, and asymptotic to, $1-2/l$ when $l$ is big.

\medskip

Recall that we assumed $\varepsilon \leq d$. Thus we get
$$
\rho(\delta, \mu) \geq
\min\Bigl( d,{l\over l+2}\Bigr)
$$

This ends the proof of the theorem.

\bigskip

\quad Let us mention a curious consequence of the previous theorem.
Let us assume that $M$ is a compact rank one surface.
This implies that there are invariant ergodic 
probability measures of full support.
Given a closed geodesic in the cylinder, we may consider the set of
ergodic measures of full support that are closest to the Dirac measure on the 
closed geodesic, with respect to the Prohorov distance.
This set is a non-empty compact set, and its elements 
should be related to the closed geodesic in some way.

\medskip

\quad
These results lead to the following problem:

\medskip 

\emph{Characterize the non-positively curved compact manifolds on which
ergodic measures of full support are dense in the set of all invariant
probability measures.}

\medskip

Recall that full support is a generic property in the set of all invariant
probability measures, as soon as the periodic geodesics are dense in $T^1M$. 
Also the set of ergodic measures is always a $G_\delta$
set. The question really is about whether ergodicity is dense or not. The
existence of embedded euclidean cylinder seems to be the only obstruction to
genericity in dimension two. In higher dimension, the question appears to be
pretty elusive.

\section{Transitivity}

We now consider the question of the transitivity of the geodesic flow
in restriction to the non-wandering set $\Omega$ of the flow. We have mentioned
that if $\Omega$ is equal to $T^1M$, then the flow is indeed transitive.
When there are wandering vectors, we have the following result.

\bigskip

{\bf Theorem 9 } \cite{cs}

\emph{Let $M$ be a rank one manifold.
We assume that there are at least three rank one periodic geodesics on $T^1M$.
Then the restriction of the geodesic flow to the closure of the set of
rank one periodic orbits is transitive.}

\bigskip

\quad 
So the density of rank one periodic vectors in the non-wandering set is
sufficient to get the transitivity of the flow. We will see below examples of
surfaces where rank one periodic orbits are not dense in $\Omega$.

\medskip 

\quad 
Here again, the requirement that there are more than two periodic orbits
rules out the case of an hyperbolic cylinder. On the unit tangent bundle of a
hyperbolic cylinder, there are exactly two periodic geodesic orbits,
corresponding to the same geometric geodesic on the surface, with its two
opposite orientations. As a consequence, the set $\Omega$ is disconnected and
transitivity does not hold on $\Omega$

\medskip 

\quad
Is the connectedness of $\Omega$ is sufficient to guarantee the
existence of a non wandering vector whose orbit is dense in $\Omega$ ? We shall
answer that question by the negative.

\medskip

\quad Let $M$ be a connected complete, non-positively curved surface, admitting
a single end isometric to a half flat cylinder
$[0,\infty[\times ({{\bf R}/{\bf Z}})$. 
Let $v$ be a periodic vector generating a geodesic whose projection
$c$ on $M$ bounds the half cylinder. We assume that the curvature is negative
outside the flat end. So the surface is separated in two parts by the curve
$c$. One part is compact and contains (the projection on $M$ of ) all rank one
periodic orbits of the geodesic flow, the other is flat and contains two
continua of rank two periodic orbits.

\medskip

\quad An example is depicted below. We take a negatively curved pant, glue two
ends together and flatten the last end so as to paste a flat cylinder. The
negatively curved part is diffeomorphic to a once punctured torus bounded by $c$.

\bigskip

\centerline{
\epsfig{figure=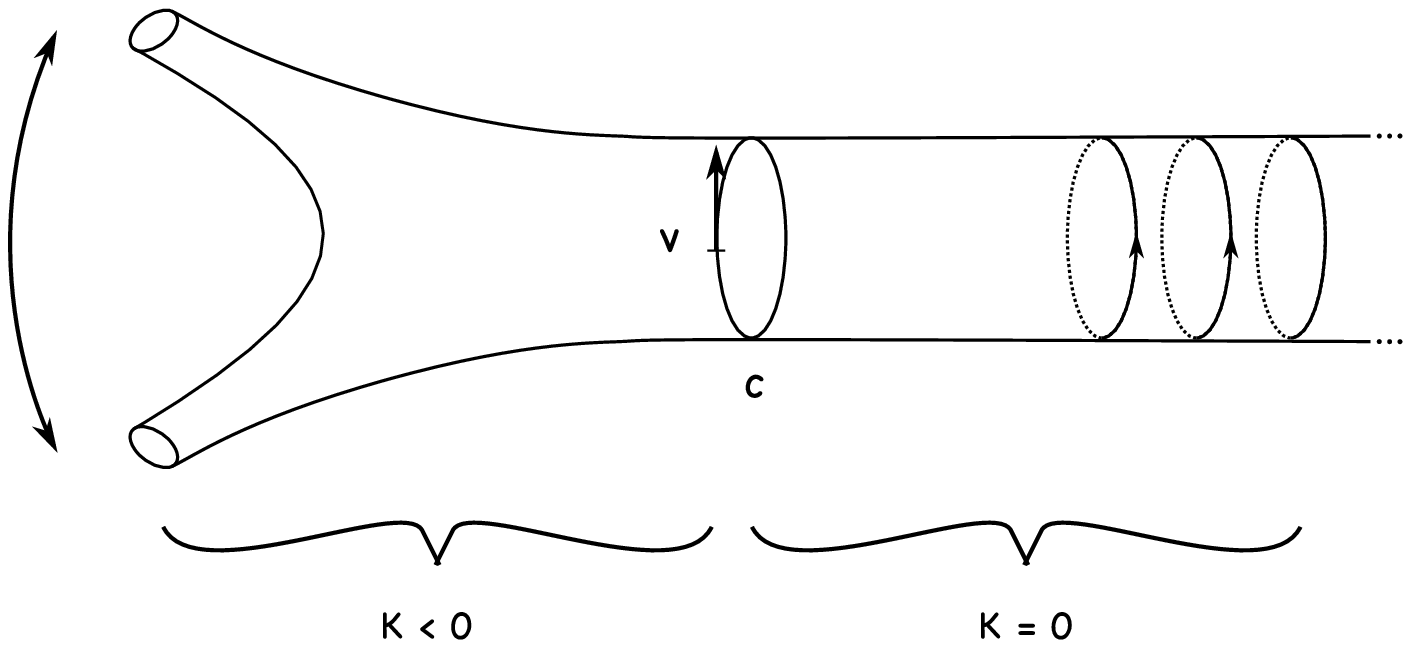,width=.8\textwidth,angle=0}}

\bigskip

{\bf Theorem 10 }

\emph{The non wandering set $\Omega\subset T^1M$ is connected, periodic orbits
  are dense in $\Omega$ and the geodesic flow is not transitive in restriction
  to $\Omega$.}

\bigskip 

{\bf Proof}\\
The cylinder contains two sets of opposite rank two periodic orbits. We will
call these orbits ``vertical'', in accordance to the figures. All other
geodesics in the cylinder go to infinity either for positive or negative time.
Any non periodic vector in the cylinder generates a geodesic that goes to
infinity either for positive or negative time. As a consequence, the only non
wandering orbits intersecting the flat cylinder are the vertical ones. See figure
below.

\bigskip

\centerline{
\epsfig{figure=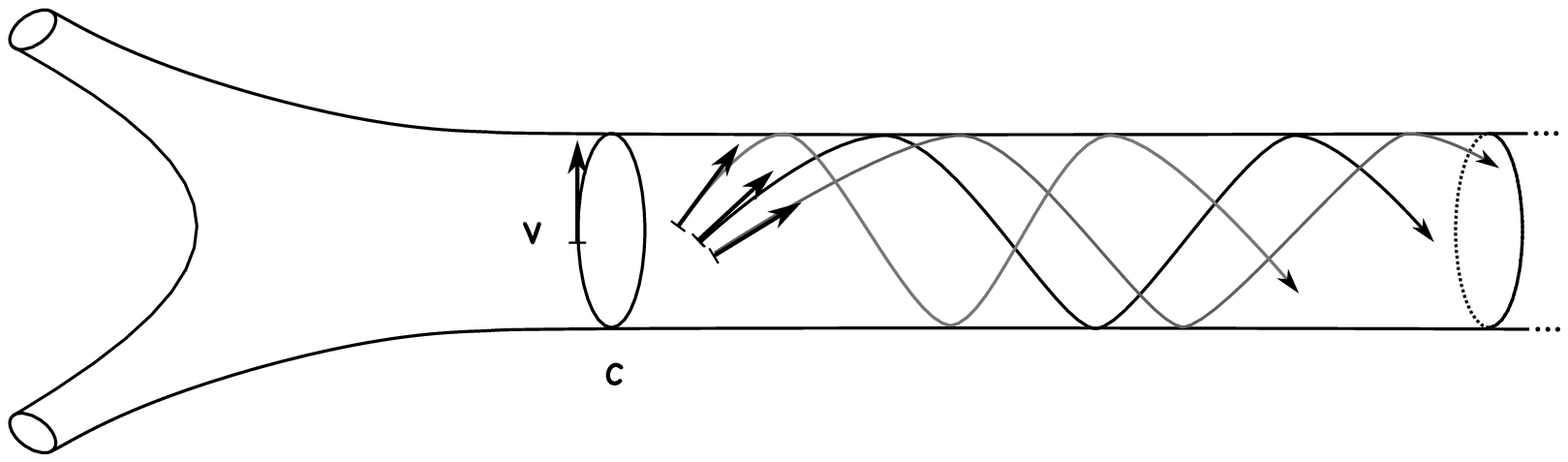,width=.8\textwidth,angle=0}}

\bigskip

\quad
This shows that rank one periodic orbits are contained in the negatively
curved part of the surface. By the rank one closing lemma, they are dense in the
set of non wandering rank one vectors. We see that the non-wandering set is
composed of the closure of the rank one periodic orbits, together with the rank
two periodic orbits contained in the cylinder.

\medskip

\quad
Also, the flow is not transitive in restriction to $\Omega$, since a non
wandering vector entering the flat cylinder is periodic, hence cannot be dense.

\medskip

\quad
It remains to show that $\Omega$ is connected. The flow is transitive in 
restriction to the closure of the rank one periodic orbits, that is, in
restriction to the subset of $\Omega$ contained in the negatively
curved part of the surface. This implies that this subset is connected.
In particular, a periodic rank one vector and its opposite are in the same
connected component of the non wandering set.

\medskip

\quad In order to prove the connectedness of $\Omega$, it is sufficient to show
that the vector $v$ on the boundary of the flat half cylinder is in the closure
of the rank one periodic vectors. Indeed, this will imply that $-v$ is also in
that closure. Now, all rank two periodic orbits in the half cylinder
can be connected by a path to $v$ or $-v$. 
So the next lemma ends the proof.

\medskip

{\bf Lemma}

\emph{The geodesic bounding the half cylinder is in the closure of the 
rank one periodic geodesics.} 

\medskip

{\bf Proof of the lemma}

The proof makes use of the action of the fundamental
group $\Gamma$ of $M$ on the boundary of its universal cover $\tilde{M}$.
Let us recall how the ideal boundary $\partial \tilde{M}$ of $\tilde{M}$ 
is defined.

\medskip 

\quad
We consider the space of half geodesics $r:{\bf R}_+\rightarrow \tilde{M}$. 
Two half geodesics $r_1$, $r_2$ are said to be asymptotic if they
stay at a bounded distance from each other: there exists $C\geq 0$ such that
$d(r_1(t), r_2(t))<C$ for all $t\geq 0$. The boundary $\partial \tilde{M}$ is
obtained by identifying asymptotic half geodesics.
We now fix some origin $x_0\in \tilde{M}$. The boundary $\partial \tilde{M}$
is homeomorphic to $T^1_{x_0}\tilde{M}$ via the map associating to each vector
of \smash{$T^1_{x_0}\tilde{M}$} the half geodesic starting from that vector.

\medskip

\quad
The limit set $\Lambda \Gamma$ is defined as the closure in 
$\partial \tilde{M}$ of the orbit of $x_0$ under $\Gamma$. It is contained 
in $\partial \tilde{M}$ and it does not depend on the choice of the origin 
$x_0$. G. Link, M. Peign\'e and J. C. Picaud showed that the end points 
of the lifts of the rank one periodic vectors are dense in
$\Lambda \Gamma \times \Lambda \Gamma$ \cite{lpp}. 
As a side note, we remark that this does not imply that the rank one periodic
orbits are dense in $\Omega$.

\medskip

\hbox{\vbox{\hsize=.5\textwidth
Let $\tilde{c}$ be a lift of the periodic geodesic $c$ that lies on the
boundary of the negatively curved part. We denote by $c^-$ and $c^+$ its two end 
points in $\Lambda\Gamma$. Let $\tilde{c}_n$ a sequence of geo\-de\-sics in
$\tilde{M}$ associated to rank one periodic geo\-de\-sics on $M$, and whose
end points $c_n^-$ and $c_n^+$ tend to $c^-$ and $c^+$. Passing to a
subsequence, we can assume that the convergence is in fact monotonous in the
neighborhood of $c^-$ and $c^+$. Let us parameterize $\tilde{c}_n$ so that
the distance from $\tilde{c}_n$ to $\tilde{c}(0)$ is realized at
$\tilde{c}_n(0)$.}
\qquad \raise 0em \vtop{
\epsfig{figure=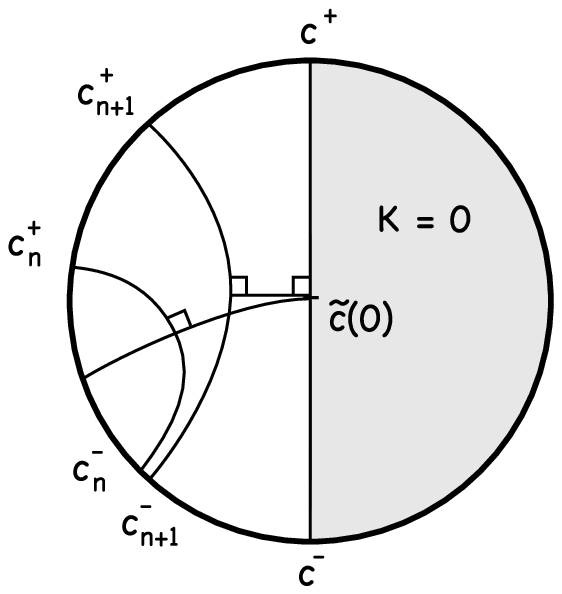,width=.4\textwidth,angle=0}
}}

\medskip

The geodesics $\tilde{c}_n$ and $\tilde{c}$ are separated by the geodesic
$\tilde{c}_{n+1}$. This implies that the sequence
$d(\tilde{c}_n(0),\tilde{c}(0))$ is decreasing. Let $w$ be an accumulation point
of the sequence $c'_n(0)$. The geodesic generated by $w$ has $c^-$ and $c^+$ as
its end points, so it bounds a flat strip together with $\tilde{c}$. Moreover,
it is in the closure of the (lifts of the) rank one periodic geodesics, so it
must be equal to $\tilde{c}$. This proves the lemma.


\end{document}